\documentclass[%
jmp,
amsmath,amssymb,
notitlepage,
superscriptaddress,
]{revtex4-1}

\usepackage{graphicx}% Include figure files
\usepackage{dcolumn}% Align table columns on decimal point
\usepackage{bm}% bold math
\usepackage[utf8]{inputenc}
\usepackage[T1]{fontenc}
\usepackage{mathptmx}

\usepackage{amsthm}
\usepackage{amsfonts}
\usepackage{mathrsfs}
\usepackage{amsmath}
\usepackage{amssymb,amsfonts}

\newtheorem{thm}{Theorem}[section]

\newtheorem{prop}{Proposition}[section]
\numberwithin{equation}{section}

\def\pf{{\textit {Proof:} }}

\newcommand{\mysection}[1]{\section{#1}\setcounter{equation}{0}}

\newfont{\bb}{msbm10 at 11pt}

\newcommand{\bal}{\begin{aligned}}      \newcommand{\eal}{\end{aligned}}
\newcommand{\ba}{\begin{array}}      \newcommand{\ea}{\end{array}}
\newcommand{\bc}{\begin{center}}     \newcommand{\ec}{\end{center}}
\newcommand{\be}{\begin{enumerate}}  \newcommand{\ee}{\end{enumerate}}
\newcommand{\beq}{\begin{eqnarray}}  \newcommand{\eeq}{\end{eqnarray}}
\newcommand{\beQ}{\begin{eqnarray*}} \newcommand{\eeQ}{\end{eqnarray*}}
\newcommand{\bi}{\begin{itemize}}    \newcommand{\ei}{\end{itemize}}
\newcommand{\bt}{\begin{tabular}}    \newcommand{\et}{\end{tabular}}
\newcommand{\bdm}{\begin{displaymath}} \newcommand{\edm}{\end{displaymath}}

\newcommand{\ls}{\setlength{\baselineskip}{12pt}
                 \setlength{\parskip}{3mm}}

\begin{document}

%\preprint{AIP/123-QED}

\title[Soliton type]{Metrics of Horowitz-Myers type with the negative constant scalar curvature}
% Force line breaks with \\

\author{Zhuobin Liang}
\thanks{\textbf{Electronic mail:} tzbliang@jnu.edu.cn}
\affiliation{Department of Mathematics, Jinan University, Guangzhou, 510632, P. R. China}

\author{Xiao Zhang}
\thanks{\textbf{Author to whom correspondence should be addressed:} xzhang@gxu.edu.cn, xzhang@amss.ac.cn}
\affiliation{ Guangxi Center for Mathematical Research, Guangxi University, Nanning, Guangxi 530004, PR China%\\This line break forced with \textbackslash\textbackslash
}%
\affiliation{Institute of Mathematics, Academy of Mathematics and Systems Science, Chinese Academy of Sciences,
Beijing 100190, PR China}%
\affiliation{School of Mathematical Sciences, University of Chinese Academy of Sciences, Beijing 100049, PR China}%

%\date{Resubmission to the revised version: October 29, 2020}

\begin{abstract}
We construct a one-parameter family of complete metrics of Horowitz-Myers type with the negative constant scalar curvature. We also verify a positive energy conjecture of
Horowitz-Myers for these metrics.
\end{abstract}

\maketitle

\mysection{Introduction}\ls

In \cite{HM}, Horowitz, Myers constructed so-called AdS solitons with toroidal topology
\begin{equation}
-r^{2}dt^{2}+\frac{1}{\frac{r^{2}}{\ell^{2}}\left(1-\frac{r_{0}^{n}}{r^{n}}\right)}dr^{2}+\frac{r^{2}}{\ell^{2}}\left(1-\frac{r_{0}^{n}}{r^{n}}\right)d\phi^{2}
+r^{2}\sum_{i=1}^{n-2}(d\theta^{i})^{2},\label{eq:ads_soliton_spacetime_metric}
\end{equation}
for some $r_{0}>0$. These metrics are globally static vacuum $(n+1)$-dimensional spacetime with cosmological constant
\beQ
\Lambda=-\frac{n(n-1)}{2\ell^{2}}.
\eeQ
The induced Riemannian metrics $g_{0}$ on the constant time slice
\begin{equation}
g_{0}=\frac{1}{\frac{r^{2}}{\ell^{2}}\left(1-\frac{r_{0}^{n}}{r^{n}}\right)}dr^{2}+\frac{r^{2}}{\ell^{2}}\left(1-\frac{r_{0}^{n}}{r^{n}}\right)d\phi^{2}+r^{2}\sum_{i=1}^{n-2}(d\theta^{i})^{2}.\label{eq:ads_soliton_t-slice_g}
\end{equation}
are ALH and are referred as Horowitz-Myers metrics. Horowitz, Myers also verified that the Hawking-Horowitz mass of $g_0$ is negative and conjectured that, among all metrics which are asymptotic to $g_0$ and with scalar curvature
\beq
S \geq -\frac{n(n-1)}{\ell^2}, \label{S}
\eeq
$g_0$ has the lowest Hawking-Horowitz mass \cite{HM}. This conjecture was shown to be true for Horowitz-Myers metrics up to order $O\big(\frac{1}{r^n}\big)$ on $R^2 \times T^{n-2}$ and with scalar curvature satisfying (\ref{S}) \cite{BCHMN}.

On the other hand, ALE or ALH metrics of Eguchi-Hanson types attracted much attentions in physics and geometry and, in particular, these metrics with constant scalar curvature play important roles in constructing examples with the negative mass \cite{EH1, EH2, L1, Z, P, L2, D, CZ}. Therefore it is natural to study the analogous question for metrics of Horowitz-Myers type. And this is the purpose of this short paper.

The paper is organized as follows. In Section II, we construct complete metrics of Horowitz-Myers type with the negative constant scalar curvature. In Section III, we show that these metrics can not develop to certain vacuum spacetimes unless the Horowitz-Myers metrics. In Section IV, we show that there is a natural almost-complex structure, which is integrable and compatible with these metrics. But the metrics are not K\"{a}hler with respect to this almost-complex structure. In Section V, we compute the Hawking-Horowitz mass and the Hamiltonian energy of these metrics and verify a positive energy conjecture of Horowitz-Myers for these metrics.

\mysection{Metrics of Horowitz-Myers type}\ls

In this section, we construct a one-parameter family complete metrics of asymptotically Horowitz-Myers type.
Let $r_{+}$ be a positive number, $V=V(r)$ be a positive smooth
function for $r>r_{+}$. For $n \geq 3$, we consider the metric
\begin{equation}
g:=\frac{1}{V}dr^{2}+Vd\phi^{2}+r^{2}\sum_{i=1}^{n-2}(d\theta^{i})^{2},\label{eq:g}
\end{equation}
where $\theta^{i}$ and $\phi$ are periodic. The periods
of $\theta^{i}$ are arbitrary, and the period of $\phi$ is
specified later. Let $\nabla$ be the Levi-Civita connection of $g$.
Denote $V^{\prime}:=\frac{dV}{dr}$. By the
curvature formulae for warped product (cf. \cite[Chapter 7]{O}), we have, for
$1\leq i,j\leq n-2$, that
\begin{equation}
\begin{aligned}\nabla_{\partial_{r}}\partial_{r} & =-\frac{1}{2}V^{-1}V^{\prime}\partial_{r},\\
\nabla_{\partial_{\phi}}\partial_{r} & =\nabla_{\partial_{r}}\partial_{\phi}=\frac{1}{2}V^{-1}V^{\prime}\partial_{\phi},\\
\nabla_{\partial_{\theta^{i}}}\partial_{r} & =\nabla_{\partial_{r}}\partial_{\theta^{i}}=\frac{1}{r}\partial_{\theta^{i}},\\
\nabla_{\partial_{\phi}}\partial_{\phi} & =-\frac{1}{2}VV^{\prime}\partial_{r},\\
\nabla_{\partial_{\theta^{i}}}\partial_{\phi} & =\nabla_{\partial_{\phi}}\partial_{\theta^{i}}=0,\\
\nabla_{\partial_{\theta^{i}}}\partial_{\theta^{j}} & =-rV\delta_{ij}\partial_{r}.
\end{aligned}
\label{eq:connection_g}
\end{equation}
Moreover, the Ricci curvature tensor $\operatorname{Ric}_{g}$ of
$g$ is diagonal with respect to the coordinate frame
\begin{equation}
\begin{aligned}\operatorname{Ric}_{g}(\partial_{r},\partial_{r}) & =-\frac{1}{2}\left(V^{\prime\prime}+\frac{n-2}{r}V^{\prime}\right)V^{-1},\\
\operatorname{Ric}_{g}(\partial_{\phi},\partial_{\phi}) & =-\frac{1}{2}\left(V^{\prime\prime}+\frac{n-2}{r}V^{\prime}\right)V,\\
\operatorname{Ric}_{g}(\partial_{\theta^{i}},\partial_{\theta^{i}}) & =-rV^{\prime}-(n-3)V
\end{aligned}
\label{eq:Ricci_g}
\end{equation}
for $1\leq i\leq n-2$, and the scalar curvature is
\[
\operatorname{Scal}_{g}=-V^{\prime\prime}-\frac{2(n-2)}{r}V^{\prime}-\frac{(n-2)(n-3)}{r^{2}}V.
\]

Now we seek the function $V$ such that $\operatorname{Scal}_{g}=S<0$. It turns out that the general solutions are
\beQ
V=\frac{r^{2}}{\ell^{2}}\left(1+\frac{a}{r^{n-1}}+\frac{b}{r^{n}}\right), \quad  \ell=\sqrt{-\frac{n(n-1)}{S}}  \label{eq:V}
\eeQ
where $a$, $b$ are two arbitrary constants. We restrict to the subset of solutions with $b=-r_{0}^{n}<0$, i.e.,
\begin{equation}
V=\frac{r^{2}}{\ell^{2}}\left(1+\frac{a}{r^{n-1}}-\frac{r_{0}^{n}}{r^{n}}\right).\label{eq:V_r0}
\end{equation}
When $a=0$, they give rise to Horowitz-Myers metrics. Since $$\lim_{r\to+\infty}V(r)=+\infty, \qquad \lim_{r\to0^{+}}V(r)=-\infty,$$
$V(r)$ possesses the largest positive root $r_{+}>0$ for any $a\in\mathbb{R}$. So $V(r)>0$ for $r>r_{+}$, and $V(r_{+})=0$. The metric $g$ has conical singularity at $r=r_{+}$. However, this singularity is removable,
provided the period of $\phi$ is specified appropriately.
\begin{thm}\label{thm1}
Let $g$ be the metric (\ref{eq:g}) with $V$ given by (\ref{eq:V_r0}) and $a\in\mathbb{R}$. Assume
that $r\geq r_{+}$, the periods of $\{\theta^{i}\}_{1\leq i\leq n-2}$ are arbitrary and the period
of $\phi$ is
\begin{equation}
\beta:=\frac{4\pi\ell^{2}}{r_{+}\left(n-1+\frac{r_{0}^{n}}{r_{+}^{n}}\right)}.\label{eq:period_phi}
\end{equation}
Then $r=r_{+}$ is a removable singularity, and $g$ is a geodesically complete metric with negative constant scalar curvature $S$ on $R^2 \times T^{n-2}$.
\end{thm}

\pf
From $V(r_{+})=0$, we have $$\frac{a}{r_{+}^{n-1}}=\frac{r_{0}^{n}}{r_{+}^{n}}-1.$$
Hence
\[
V^{\prime}(r_{+})=\frac{r_{+}}{\ell^{2}}\left[2-(n-3)\frac{a}{r_{+}^{n-1}}-(n-2)\frac{b}{r_{+}^{n}}\right]=\frac{4\pi}{\beta}>0.
\]
By Implicit Function Theorem, $V$ has smooth inverse
$\tilde{V}$ around $r_{+}$. Set $\rho=V^{\frac{1}{2}}(r)$ for $r>r_{+}$.
Then $$r=\tilde{V}(\rho^{2}), \qquad d\rho=\frac{1}{2}V^{-\frac{1}{2}}V^{\prime}dr.$$
Replacing $r$ with $\rho$, we rewrite $g$ as follows
\[
g=\frac{4}{V^{\prime}(r_{+})^{2}} \left[h(\rho^{2})\rho^{2}d\rho^{2}+\left(d\rho^{2}+\rho^{2}d\Phi^{2}\right) \right] +\tilde{V}^2(\rho^{2})\sum_{i=1}^{n-2}(d\theta^{i})^{2},
\]
where $\Phi=\frac{2\pi\phi}{\beta}$ and
\[
h(s):=s^{-1}\left[\frac{V^{\prime}(r_{+})^{2}}{\left(V^{\prime}\circ \tilde{V}(s)\right)^{2}}-1\right].
\]
Clearly, $\Phi$ is of period $2\pi$, and $h(s)$ is smooth around
$s=0$. Now we apply coordinate transformation $x=\rho\cos\Phi$ and
$y=\rho\sin\Phi$, and find that $g$ becomes
\[
\frac{4}{V^{\prime}(r_{+})^{2}}\left[h(x^{2}+y^{2})(xdx+ydy)^{2}+\left(dx^{2}+dy^{2}\right)\right]+\tilde{V}^{2}(x^{2}+y^{2})\sum_{i=1}^{n-2}(d\theta^{i})^{2}.
\]
From this, we see that $g$ is smooth at $(x,y)=(0,0)$ and $r=r_{+}$ is a removable singularity. Thus $g$ is
a complete metric and the underlying manifold has the topology $R ^2\times T^{n-2}$. Moreover, $g$
has constant scalar curvature $S$. So the proof is complete.

Consequently, the initial data $(g,K\equiv0)$ satisfying the constraint
equations for the vacuum Einstein equations with negative cosmological
constant $\Lambda=\frac{1}{2}S$.

\section{Uniqueness}

In \cite{GSW}, Galloway, Surya and Woolgar established a uniqueness theorem
of AdS solitons for static metrics which can be conformally compactified,
see also \cite{ACD} for certain extension. In their case some asymptotic
behavior of the lapse $N$ at spatial infinity is required. In the following
we prove that the metric $g$ constructed in Theorem \ref{thm1} can not develop
to certain vacuum spacetime unless $a=0$, where the lapse $N$ does not require to
satisfy any asymptotic condition at spatial infinity.

\begin{thm}\label{thm2}
Let $g$ be the metric constructed in Theorem \ref{thm1}. Let $N=N(r,\phi,\theta^{i})$
be a positive function where $r>r_{+}$. Suppose the metric
\[
\widetilde{g}=-N^{2}dt^{2}+\frac{1}{V}dr^{2}+Vd\phi^{2}+r^{2}\sum_{i=1}^{n-2}(d\theta^{i})^{2}.
\]
satisfies the vacuum Einstein field equations
\begin{equation}
\operatorname{Ric}_{\widetilde{g}}-\frac{1}{2}\operatorname{Scal}_{\widetilde{g}}\widetilde{g}+\Lambda\widetilde{g}=0\label{eq:VEE_Lambda}
\end{equation}
for the negative cosmological constant $\Lambda<0$. Then it holds that
\beQ
N=c r, \quad \ell^{2}=-\frac{n(n-1)}{2\Lambda}, \quad V=\frac{r^{2}}{\ell^{2}}\left(1-\frac{r_{0}^{n}}{r^{n}}\right)
\eeQ
where $c$ is some positive constant which can be chosen as $1$ by defining new time $\tilde{t}=c t$. Therefore $\widetilde{g}$ is a AdS soliton.
\end{thm}

\pf
Let $g$ be the spatial part of $\widetilde{g}$. By (\ref{eq:VEE_Lambda}),
we have
\begin{equation}
\nabla^{2}N=N\left(\operatorname{Ric}_{g}-\frac{2\Lambda}{n-1}g\right).\label{eq:Hessian_N}
\end{equation}

\textit{Claim 1: There exist smooth functions
\[
k=k(r), \quad h=h(\phi), \quad f_{i}=f_{i}(\theta^{i})
\]
where $h$ and $f_{i}$ are periodic, $1\leq i\leq n-2$, such that
\begin{equation}
N=k(r)+V^{\frac{1}{2}}(r)h(\phi)+r\sum_{i=1}^{n-2}f_{i}(\theta^{i}).\label{eq:N_decomposition}
\end{equation}
Proof of the Claim 1}: Recall that $\operatorname{Ric}_{g}$
is diagonal. It follows from (\ref{eq:Hessian_N}) that
$\nabla^{2}N$ is also diagonal.
From this and (\ref{eq:connection_g}), we have
\beQ
\nabla^{2}N(\partial_{r},\partial_{\phi})=\nabla^{2}N(\partial_{r},\partial_{\theta^{i}})=
\nabla^{2}N(\partial_{\phi},\partial_{\theta^{i}})=
\nabla^{2}N(\partial_{\theta^{i}},\partial_{\theta^{j}})=0
\eeQ
for $1\leq i,j\leq n-2$ ($i\neq j$), we obtain
\beQ
\partial_{r}\partial_{\phi}(V^{-\frac{1}{2}}N)=
\partial_{r}\partial_{\theta^{i}}(r^{-1}N)=
\partial_{\phi}\partial_{\theta^{i}}N=
\partial_{\theta^{i}}\partial_{\theta^{j}}N=0.
\eeQ
Therefore $N$ must take
the form (\ref{eq:N_decomposition}). Since $\phi$ and $\theta^{i}$
are periodic, $h$ and $f_{i}$ are periodic.

\textit{Claim 2: $N$ depends only on $r$.\\\\
Proof of the Claim 2}: By (\ref{eq:Hessian_N}), we have
\[
\begin{aligned}
\partial_{\phi}\partial_{\phi}N-\nabla_{\partial_{\phi}}\partial_{\phi}(N)&=\nabla^{2}N(\partial_{\phi},\partial_{\phi})=N\left(\operatorname{Ric}_{g}(\partial_{\phi},\partial_{\phi})-\frac{2\Lambda}{n-1}g(\partial_{\phi},\partial_{\phi})\right).
\end{aligned}
\]
By (\ref{eq:connection_g}) and (\ref{eq:Ricci_g}), we obtain
\begin{equation}
\partial_{\phi}\partial_{\phi}N+\frac{1}{2}VV^{\prime}\partial_{r}N=-N\left(\frac{1}{2}V^{\prime\prime}+\frac{n-2}{2r}V^{\prime}+\frac{2\Lambda}{n-1}\right)V.\label{eq:pf_static_2}
\end{equation}
Taking the derivative $\partial_{\phi}$ on both sides, and noting that $\partial_{\phi}N=V^{\frac{1}{2}}\frac{dh}{d\phi}$ by (\ref{eq:N_decomposition}), we then obtain,
\begin{equation}
-4\frac{d^{3}h}{d\phi^{3}}=\frac{dh}{d\phi}\cdot\left[2\left(V^{\prime\prime}+\frac{n-2}{r}V^{\prime}\right)V+(V^{\prime})^{2}+\frac{8\Lambda}{n-1}V\right].\label{eq:h_derivative}
\end{equation}
Note the term in the bracket is equal to
\beq
\begin{aligned}
\frac{r^{2}}{\ell^{4}}\left[4\left(n+\frac{2\Lambda\ell^{2}}{n-1}\right) +(n-3)^{2} \frac{a^{2}}{r^{2(n-1)}}+n(n-2)\frac{b^{2}}{r^{2n}}
+2(n-2)^{2}\frac{ab}{r^{2n-1}}+8\left(\frac{\Lambda\ell^{2}}{n-1}+1\right)\frac{a}{r^{n-1}}
+\left(2n+\frac{8\Lambda\ell^{2}}{n-1}\right)\frac{b}{r^{n}}\right]   \label{eq:h_derivative_2}
\end{aligned}
\eeq
where $b=-r_{0}^{n}$. Now we assert that $\frac{dh}{d\phi}\equiv0$. If $\frac{dh}{d\phi}(\phi_{0})\neq0$ for some $\phi_{0}$. As $h$ depends only on $\phi$, we conclude that (\ref{eq:h_derivative_2}) must be constant
for $r>r_{+}$. This implies that $$n+\frac{2\Lambda\ell^{2}}{n-1}=b=0.$$
But $b=-r_{0}^{n}\neq0$, this yields a contradiction. Consequently, $\frac{dh}{d\phi}\equiv0$, and $\partial_{\phi}N=0$.

Similarly, we can show that $\partial_{\theta^{i}}N=0$. Indeed, by (\ref{eq:Hessian_N}), we obtain
\beQ
\begin{aligned}
\partial_{\theta^{i}}\partial_{\theta^{i}}N-\nabla_{\partial_{\theta^{i}}}\partial_{\theta^{i}}(N)&=\nabla^{2}N(\partial_{\theta^{i}},\partial_{\theta^{i}})
=N\left(\operatorname{Ric}_{g}(\partial_{\theta^{i}},\partial_{\theta^{i}})-\frac{2\Lambda}{n-1}g(\partial_{\theta^{i}},\partial_{\theta^{i}})\right)
\end{aligned}
\eeQ
Using (\ref{eq:connection_g}) and (\ref{eq:Ricci_g}), we obtain
\begin{equation}
\partial_{\theta^{i}}\partial_{\theta^{i}}N+rV\partial_{r}N=-Nr^{2}\left(\frac{1}{r}V^{\prime}+\frac{n-3}{r^{2}}V+\frac{2\Lambda}{n-1}\right).\label{eq:pf_static_3}
\end{equation}
Taking the derivative $\partial_{\theta^{i}}$ on both sides, and noting that (\ref{eq:N_decomposition}) gives $\partial_{\theta^{i}}N=r\frac{df_{i}}{d\theta^{i}}$, we get,
\beq
\frac{d^{3}f_{i}}{d(\theta^{i})^{3}}=-\frac{df_{i}}{d\theta^{i}}\cdot r^{2}\left(\frac{1}{r}V^{\prime}+\frac{n-2}{r^{2}}V+\frac{2\Lambda}{n-1}\right)
=-\frac{df_{i}}{d\theta^{i}}\cdot\frac{r^{2}}{\ell^{2}}\cdot\left(n+\frac{2\Lambda\ell^{2}}{n-1}+\frac{a}{r^{n-1}}\right).
\label{eq:fi_derivative-1}
\eeq
Now we assert that $\frac{df_{i}}{d\theta^{i}}\equiv0$. If it is not true, same as before, we can conclude that
\begin{equation}
\frac{r^{2}}{\ell^{2}}\cdot\left(n+\frac{2\Lambda\ell^{2}}{n-1}+\frac{a}{r^{n-1}}\right)=C, \quad r>r_{+}.\label{eq:pf_static_1}
\end{equation}
for certain constant $C$.

Case 1: $n>3$. Note that (\ref{eq:pf_static_1}) gives that
\[
n+\frac{2\Lambda\ell^{2}}{n-1}=a=0.
\]
Then (\ref{eq:fi_derivative-1}) becomes $\frac{d^{3}f_{i}}{d(\theta^{i})^{3}}\equiv0$.
Since $f_{i}$ is periodic, this implies $f_{i}$ must be a constant,
which yields a contradiction as $\frac{df_{i}}{d\theta^{i}}(\theta_{0}^{i})\neq0$.
Consequently, $\frac{df_{i}}{d\theta^{i}}\equiv0$, and
$\partial_{\theta^{i}}N=0$.

Case 2: $n=3$. In this case, $\frac{df_{1}}{d\theta^{1}}(\theta_{0}^{1})\neq0$, and (\ref{eq:pf_static_1}) only gives that
\[
n+\frac{2\Lambda\ell^{2}}{n-1}=0\quad \Rightarrow \quad \Lambda = -\frac{3}{\ell^2}.
\]
Since $\partial_{\phi}N=0$, (\ref{eq:pf_static_2}) gives $\partial_{r}N=\frac{1}{r}N$. When $n=3$, (\ref{eq:N_decomposition}) implies $\frac{dh}{d\phi}=0$, we obtain $N=\kappa r+rf_{1}(\theta^{1})$
for some constant $\kappa$. Putting this into (\ref{eq:pf_static_3}) with $i=1$, we obtain
\begin{equation}
\frac{d^{2}f_{1}}{d(\theta^{1})^{2}}+\frac{a}{\ell^{2}}f_{1}=-\frac{a\kappa}{\ell^{2}}.\label{eq:pf_static_4}
\end{equation}

If $a\equiv0$, then (\ref{eq:pf_static_4}) becomes $\frac{d^{2}f_{1}}{d(\theta^{1})^{2}}\equiv0$.
Since $f_{1}$ is periodic, this implies $f_{1}$ must be a constant,
which yields a contradiction as $\frac{df_{1}}{d\theta^{1}}(\theta_{0}^{1})\neq0$.

If $a<0$, then the solutions of (\ref{eq:pf_static_4}) are
\[
f_{1}(\theta^{1})=-\kappa+A_{1}e^{\frac{\sqrt{-a}}{\ell}\theta^{1}}+A_{2}e^{-\frac{\sqrt{-a}}{\ell}\theta^{1}},
\]
with $A_{1},A_{2}\in\mathbb{R}$. Since $f_{1}$ is periodic, this
implies $A_{1}=A_{2}=0$ and so $f_{1}=-\kappa$. Again, this is
impossible as $\frac{df_{1}}{d\theta^{1}}(\theta_{0}^{1})\neq0$.

If $a>0$, then the solutions of (\ref{eq:pf_static_4}) are
\[
f_{1}(\theta^{1})=-\kappa+A_{1}\sin\left(\frac{\sqrt{-a}}{\ell}\theta^{1}+A_{2}\right)
\]
with $A_{1},A_{2}\in\mathbb{R}$. It follows
\[
N=A_{1}r\sin\left(\frac{\sqrt{-a}}{\ell}\theta^{1}+A_{2}\right).
\]
But this is impossible since $N>0$ for all $\theta^{1}$.

As a result, we conclude that $\frac{df_{1}}{d\theta^{1}}\equiv0$, so $\partial_{\theta^{1}}N=0$. This proves
the claim.

Now we obtain
\begin{equation}
\left\{ \begin{aligned}\operatorname{Ric}_{\widetilde{g}}(\partial_{t},\partial_{t}) & =NN^{\prime\prime}V+NN^{\prime}\left(V^{\prime}+\frac{n-2}{r}V\right),\\
\operatorname{Ric}_{\widetilde{g}}(\partial_{r},\partial_{r}) & =-N^{-1}N^{\prime\prime}-\frac{1}{2}N^{-1}N^{\prime}V^{-1}V^{\prime}-\frac{1}{2}V^{-1}\left(V^{\prime\prime}+\frac{n-2}{r}V^{\prime}\right),\\
\operatorname{Ric}_{\widetilde{g}}(\partial_{\phi},\partial_{\phi}) & =-\frac{1}{2}N^{-1}N^{\prime}VV^{\prime}-\frac{1}{2}V\left(V^{\prime\prime}+\frac{n-2}{r}V^{\prime}\right),\\
\operatorname{Ric}_{\widetilde{g}}(\partial_{\theta^{i}},\partial_{\theta^{i}}) & =-rV^{\prime}-(n-3)V-N^{-1}N^{\prime}rV,
\end{aligned}
\right.\label{eq:Ric_spacetime}
\end{equation}
for $1\leq i\leq n-2$. Since Einstein equations (\ref{eq:VEE_Lambda}) gives that
\[
\operatorname{Ric}_{\widetilde{g}}(V^{\frac{1}{2}}\partial_{r},V^{\frac{1}{2}}\partial_{r})-\operatorname{Ric}_{\widetilde{g}}(V^{-\frac{1}{2}}\partial_{\phi},V^{-\frac{1}{2}}\partial_{\phi})=0.
\]
Substituting this into (\ref{eq:Ric_spacetime}), we get $N^{-1}N^{\prime\prime}V=0$. This implies $N^{\prime\prime}=0$, so $N=cr+d$ for some $c,d\in\mathbb{R}$. Putting this into the following equation
\[
\operatorname{Ric}_{\widetilde{g}}(\partial_{t},\partial_{t})=\frac{2\Lambda}{n-1}\widetilde{g}(\partial_{t},\partial_{t})=-\frac{2\Lambda}{n-1}N^{2},
\]
and applying the explicit expression for $\operatorname{Ric}_{\widetilde{g}}(\partial_{t},\partial_{t})$
in (\ref{eq:Ric_spacetime}), we have
\[
-\frac{c}{cr+d}\left(V^{\prime}+\frac{n-2}{r}V\right)=\frac{2\Lambda}{n-1}, \quad
V^{\prime}+\frac{n-2}{r}V=\frac{r}{\ell^{2}}\left(n+\frac{a}{r^{n-1}}\right).
\]
Therefore
\[
-\frac{cr}{cr+d}\left(n+\frac{a}{r^{n-1}}\right)=\frac{2\Lambda\ell^{2}}{n-1}.
\]
for all $r>r_{+}$. It follows $$a=d=0, \quad \Lambda=-\frac{n(n-1)}{2\ell^{2}}, \quad N=cr,\quad c>0.$$
Thus the proof is completed.

\section{Complex structure}

The natural complex structures for metrics of Eguchi-Hanson types were studied in \cite{P, L2, CZ}.
In this section, we assume dimension $n=2+2k$ for some $k\geq 1$.
We show that, on the metric $g$ constructed in Theorem \ref{thm1}, there exists a natural
almost-complex structure $J$, which is integrable and $g$-compatible.
However $g$ is not K\"{a}hler with respect to this $J$.

For $r>r_{+}$, we define the almost-complex structure $J$ by
\begin{equation}
\frac{1}{V}dr\longmapsto d\phi,\quad d\theta^{j}\longmapsto d\theta^{k+j}, \quad  1\leq j\leq k. \label{eq:J_def}
\end{equation}

\begin{prop}
It holds that $J$ extends smoothly to $r=r_{+}$. So $J$ is an almost-complex
structure defined on the entire initial data.
\end{prop}

\pf
Recall the coordinates $(\rho,\Phi)$ and $(x,y)$ which are introduced
in the proof of Theorem \ref{thm1}. We rewrite
\begin{align*}
\frac{1}{V}dr & =\frac{2}{V^{\prime}\circ V^{-1}(\rho^{2})}\cdot\frac{1}{\rho}d\rho =\frac{2}{V^{\prime}\circ V^{-1}(\rho^{2})}\cdot\frac{1}{\rho^{2}}(xdx+ydy)
\end{align*}
and
\[
d\phi=\frac{\beta}{2\pi}d\Phi=\frac{\beta}{2\pi}\cdot\frac{1}{\rho^{2}}(-ydx+xdy).
\]
Note that $\beta$ is defined by (\ref{eq:period_phi}). So we have
\[
J(xdx+ydy)=\frac{\beta}{2\pi}\cdot\frac{V^{\prime}\circ V^{-1}(\rho^{2})}{2}(-ydx+xdy).
\]
Letting
\[
u(s)=\frac{\beta}{2\pi}\cdot\frac{V^{\prime}\circ V^{-1}(s)}{2},
\]
we rewrite
\[
J(xdx+ydy)=u(\rho^{2})(-ydx+xdy).
\]
By the proof of Theorem \ref{thm1}, it is easy to see that $u(s)$
is smooth positive function around $s=0$ and it satisfies $u(0)=1$.
From the above equation and from the fact that $J^{2}=-\operatorname{Id}$,
we see that $J$ maps $\text{span}\{dx,dy\}$ onto $\text{span}\{dx,dy\}$.
Denoting by $A$ the transformation matrix so that
\[
J(dx,dy)=(dx,dy)A,
\]
we find
\[
\begin{aligned}A= & \frac{1}{\rho^{2}}\left(\begin{array}{cc}
xy\big(\frac{1}{u(\rho^{2})}-u(\rho^{2})\big) & \big(1-\frac{1}{u(\rho^{2})}\big)x^{2}+\big(1-u(\rho^{2})\big)y^{2}\\
\big(u(\rho^{2})-1\big)x^{2}+\big(\frac{1}{u(\rho^{2})}-1\big)y^{2} & xy\big(u(\rho^{2})-\frac{1}{u(\rho^{2})}\big)
\end{array}\right)
+\left(\begin{array}{cc}
 & -1\\
1
\end{array}\right).
\end{aligned}
\]
Since $u(s)$ is smooth around $s=0$ and $u(0)=1$, it is easy to
see that $A$ extends smoothly to $\rho=0$ and it holds $A=\left(\begin{array}{cc}
 & -1\\
1
\end{array}\right)$ at $\rho=0$. So the proof is complete.

\begin{prop}
The almost-complex structure $J$ is integrable and $g$-compatible.
\end{prop}

\pf
With respect to $J$ given by (\ref{eq:J_def}), the $(1,0)$-forms
are spanned by
\[
\frac{1}{V}dr+id\phi,\quad d\theta^{j}+id\theta^{k+j}, \quad 1\leq j\leq k.
\]
It is easy to see that all of these are closed forms. So $J$ is integrable. Now we compute
\begin{align*}
g(J\cdot,J\cdot) & =\frac{1}{V}J(dr)\otimes J(dr)+VJ(d\phi)\otimes J(d\phi)+r^{2}\sum_{j=1}^{2k}J(d\theta^{j})\otimes J(d\theta^{j})\\
 & =Vd\phi\otimes d\phi+\frac{1}{V}dr\otimes dr+r^{2}\sum_{j=1}^{2k}d\theta^{j}\otimes d\theta^{j} =g.
\end{align*}
Therefore $J$ is $g$-compatible and the proof is complete.

Since the form
\begin{align*}
g(\cdot,J\cdot) & =\frac{1}{V}dr\otimes J(dr)+Vd\phi\otimes J(d\phi)+r^{2}\sum_{j=1}^{2k}d\theta^{j}\otimes J(d\theta^{j})\\
 & =dr\otimes d\phi-d\phi\otimes dr+r^{2} \sum_{j=1}^{k}\left(d\theta^{j}\otimes d\theta^{k+j}-d\theta^{k+j}\otimes d\theta^{j}\right)\\
 & =dr\wedge d\phi+r^{2}\sum_{j=1}^{k}d\theta^{j}\wedge d\theta^{k+j}
\end{align*}
is not closed, $g$ is not K\"{a}hler with respect to $J$.

\mysection{Energy}

There are several definitions of the total energy for asymptotically Horowitz-Myers metrics, e.g., the Hawking-Horowitz mass used in \cite{HM},
and the Ashtekar-Magnon mass used in \cite{GSW}. Two masses are the same for Horowitz-Myers metrics \cite{GSW}.
In this section, we compute the total energy of the metric $g$ constructed in Theorem \ref{thm1}. We show that the Hamiltonian
energy defined in \cite{BCHMN} is equal to the Hawking-Horowitz mass up to certain constant.
We also verify a positive energy conjecture of Horowitz-Mayers for metric $g$.

Firstly, we compute the Hawking-Horowitz mass of $g$. Let $T_{r}^{n-1}$ be the constant $r$ slice in metric $g$.
Recall that the period of $\phi$ is $\beta$. The period of $\theta^{i}$ is denoted as $\lambda_{i}$, $1\leq i\leq n-2$.
Let $H$ be the mean curvature of $T_{r}^{n-1}$ with respect to the unit normal $V^{\frac{1}{2}}\partial_{r}$. We obtain
\begin{align*}
H & =\operatorname{div}_{T_{r}^{n-1}}(V^{\frac{1}{2}}\partial_{r})\\
 & =\big<\nabla_{V^{-\frac{1}{2}}\partial_{\phi}}(V^{\frac{1}{2}}\partial_{r}),V^{-\frac{1}{2}}\partial_{\phi}\big>+\sum_{i=1}^{n-2}\big<\nabla_{r^{-1}\partial_{\theta^{i}}}(V^{\frac{1}{2}}\partial_{r}),r^{-1}\partial_{\theta^{i}}\big>\\
 & =\frac{1}{2}V^{-\frac{1}{2}}\partial_{r}g_{\phi\phi}+\frac{1}{2}\sum_{i=1}^{n-2}r^{-2}V^{\frac{1}{2}}\partial_{r}g_{\theta^{i}\theta^{i}}\\
 & =V^{-\frac{1}{2}}\left[\frac{1}{2}V^{\prime}+(n-2)r^{-1}V\right].
\end{align*}
Since $V$ is given by (\ref{eq:V_r0}), we further obtain
\begin{equation}
H=\frac{n-1}{\ell}+\frac{r_{0}^{n}}{2\ell r^{n}}+O(r^{-2(n-1)}).\label{eq:mean_curvature_r_constant}
\end{equation}

The reference spacetime metric is chosen as
\begin{equation}
-r^{2}dt^{2}+\frac{\ell^{2}}{r^{2}}dr^{2}+\frac{r^{2}}{\ell^{2}}d\phi^{2}+r^{2}\sum_{i=1}^{n-2}(d\theta^{i})^{2}\label{eq:ads-g}
\end{equation}
by taking $r_{0}=0$ in (\ref{eq:ads_soliton_spacetime_metric}), where the period of $\theta^{i}$ is $\lambda_{i}$ and the period of $\phi$ is $2\pi \ell^2$. Denote
\begin{eqnarray}
\lambda:=\lambda_{1}\cdot\dots\cdot\lambda_{n-2}
\end{eqnarray}
the volume of the torus $T^{n-2} _{\theta ^i}$. The perimeter of $T_\phi$ with respect to the metric $\frac{1}{\ell^2} d\phi ^2$ is
$2\pi \ell$. Let
\begin{eqnarray}
\breve{g}=\frac{\ell^{2}}{r^{2}}dr^{2}+\frac{r^{2}}{\ell^{2}}d\phi^{2}+r^{2}\sum_{i=1}^{n-2}(d\theta^{i})^{2}\label{breve-g}
\end{eqnarray}
be the induced Rimannian metric of constant time slice in (\ref{eq:ads-g}). Note that $\breve{g}$
is of constant sectional curvature $-1/\ell^{2}$. Let $H_{0}$ be the
mean curvature of $T_{r}^{n-1}$ with respect to $\breve{g}$. Letting
$a=0$ and $r_{0}=0$ in (\ref{eq:mean_curvature_r_constant}), we
obtain
\begin{equation}
H_{0}=\frac{n-1}{\ell}.\label{eq:mean_curvature_r_constant_ref}
\end{equation}

Recall the lapse of reference metric (\ref{eq:ads-g})
\begin{eqnarray}
N=r.
\end{eqnarray}
The Hakwing-Horowitz mass of $g$ is then defined as \cite{HM}
\begin{equation}
E_{{\tiny\mbox{HH}}}(g)=-\frac{1}{8\pi G}\lim_{r\to\infty}\int_{T_{r}^{n-1}}N(H-H_{0})dv_{g_{|T_{r}^{n-1}}}.\label{eq:HM-energy}
\end{equation}
where $G$ is the $(n+1)$-dimensional Newton's constant. Note that the area of $T^{n-1}_{r}$ is $\lambda\beta r^{n-2}V^{\frac{1}{2}}$.
Inserting this and equations (\ref{eq:mean_curvature_r_constant}) and (\ref{eq:mean_curvature_r_constant_ref})
into (\ref{eq:HM-energy}), we obtain
\begin{equation}
\begin{aligned}E_{{\tiny\mbox{HH}}}(g) & =-\frac{1}{8\pi G}\lim_{r\to\infty}\lambda\beta V^{\frac{1}{2}}r^{n-1}\left(\frac{r_{0}^{n}}{2\ell r^{n}}+O(r^{-2(n-1)})\right) =-\frac{\lambda\beta r_{0}^{n}}{16\pi G\ell^{2}}.
\end{aligned}
\label{eq:HM-energy_1}
\end{equation}
Since $\beta$ is given by (\ref{eq:period_phi}), we obtain
\begin{equation}
E_{{\tiny\mbox{HH}}}(g)=-\frac{\lambda r_{0}^{n}}{4Gr_{+}\left(n-1+\frac{r_{0}^{n}}{r_{+}^{n}}\right)}.\label{eq:HM-energy_2}
\end{equation}
From (\ref{eq:HM-energy_1}) or (\ref{eq:HM-energy_2}), one sees immediately that $$E_{{\tiny\mbox{HH}}}(g)<0$$ for all $a\in\mathbb{R}$.

Now we provide the Hamiltonian energy of an asymptotically Horowitz-Myers metric $g$ given by Barzegar, Chru\'{s}ciel, H\"{o}rzinger, Maliborski and Nguyen \cite{BCHMN}. Let
\[
\breve{e}_{1}=\frac{r}{\ell}\partial_{r},\quad\breve{e}_{2}=\frac{\ell}{r}\partial_{\phi},
\quad\breve{e}_{i}=\frac{1}{r}\partial_{\theta^{i-2}},\quad i=3,\dots,n.
\]
Then $\{\breve{e}_{i}\}_{i=1}^{n}$ is an orthonormal frame of $\breve{g}$. Let $\{\breve{e} ^i \} _{i=1} ^{n}$ be its dual frame.
Let $\breve{\nabla}$ be the Levi-Civita connection of $\breve{g}$. It is easy to compute
\[
\breve{\nabla}_{\breve{e}_{1}}\breve{e}_{1}=0,\quad
\breve{\nabla}_{\breve{e}_{2}}\breve{e}_{1}=\frac{1}{r}\partial_{\phi},\quad
\breve{\nabla}_{\breve{e}_{2}}\breve{e}_{2}=-\frac{r}{\ell^{2}}\partial_{r}
\]
and, for $3\leq i\leq n$,
\[
\breve{\nabla}_{\breve{e}_{i}}\breve{e}_{1}=\frac{1}{r\ell}\partial_{\theta^{i-2}},\quad
\breve{\nabla}_{\breve{e}_{i}}\breve{e}_{i}=-\frac{r}{\ell^{2}}\partial_{r}.
\]
These identities follow from (\ref{eq:connection_g}) with $V=\frac{r^{2}}{\ell^{2}}$.

Denote $g_{ij}=g(\breve{e}_{i},\breve{e}_{j})$ and $a_{ij}=g_{ij}-\delta_{ij}$, $1\leq i,j\leq n$, on each end. Suppose that $a_{ij}$ satisfy the fall-off conditions
\beq
a_{ij}=O\big(\frac{1}{r^n}\big), \quad \breve{\nabla} _k a_{ij} =O\big(\frac{1}{r^n}\big), \quad \breve{\nabla}_l \breve{\nabla} _k a_{ij} =O\big(\frac{1}{r^n}\big),\label{aij}
\eeq
and the scalar curvature $\operatorname{Scal}_g$ satisfies
\beq
\int_M \Big(\operatorname{Scal}_g+\frac{n(n-1)}{\ell^2}\Big) N \breve{e}^1 \wedge \breve{e}^2 \wedge \cdots \wedge \breve{e}^n < \infty.
\eeq
For this metric, the Hamiltonian energy up to certain constant is \cite{BCHMN}
\begin{equation}
E(g)=\frac{1}{4\mbox{Vol}(\breve{T}^{n-1})}\lim _{r \rightarrow \infty}\int _{T_r ^{n-1} }\mathcal{E} N \breve{e}^2 \wedge \cdots \wedge \breve{e}^n, \label{E2}
\end{equation}
where $\mbox{Vol}(\breve{T}^{n-1})=2\pi \ell \lambda$, and
\[
\mathcal{E}=\breve{\nabla}^{i}g_{1i}-\breve{\nabla}_{1}\operatorname{tr}_{\breve{g}}(g)-\frac{1}{\ell}(a_{11}-g_{11}\operatorname{tr}_{\breve{g}}(a)).
\]

Under the frame $\{\breve{e}_{i}\}_{i=1}^{n}$, the metric $g$ constructed in Theorem \ref{thm1} satisfies
\beq
\begin{aligned}
(g_{ij})&=\text{diag}\left(\left(1+\frac{a}{r^{n-1}}-\frac{r_{0}^{n}}{r^{n}}\right)^{-1}, \left(1+\frac{a}{r^{n-1}}-\frac{r_{0}^{n}}{r^{n}}\right),1,\dots,1\right),\\
(a_{ij})&=\text{diag}\left(-\frac{a}{r^{n-1}}+\frac{r_0 ^n}{r^n}+o(r^{n+1}), \frac{a}{r^{n-1}}-\frac{r_0 ^n}{r^n}, 0, \dots,0 \right)\label{aij-g}
\end{aligned}
\eeq
Therefore, $g$ does not satisfy the condition (\ref{aij}). But its Hamiltonian energy is still finite. Indeed, we can prove
\begin{prop}\label{E-equal}
For metrics of Horowitz-Myers type constructed in Theorem \ref{thm1},
\beQ
E_{\tiny\mbox{HH}}(g)=\frac{\lambda\ell}{2G}E(g).
\eeQ
\end{prop}

\pf By a direct computation, we have
\begin{align*}
\breve{\nabla}^{i}g_{1i} =(n+1)\frac{r^{2}}{\ell^{3}}V^{-1}-\frac{r^{3}}{\ell^{3}}V^{-2}V^{\prime}-\frac{\ell}{r^{2}}V-\frac{n-2}{\ell} =-\frac{a}{\ell r^{n-1}}+O(r^{-(n+1)}),
\end{align*}
\begin{align*}
\breve{\nabla}_{1}\operatorname{tr}_{\breve{g}}(g)  =\frac{r}{\ell}\partial_{r}\left(\frac{r^{2}}{\ell^{2}}V^{-1}+\frac{\ell^{2}}{r^{2}}V+n-2\right)=O(r^{-(n+1)}),
\end{align*}
and
\begin{align*}
\frac{1}{\ell}(a_{11}-g_{11}\operatorname{tr}_{\breve{g}}(a)) =\frac{3r^{2}}{\ell^{3}}V^{-1}-\frac{r^{4}}{\ell^{5}}V^{-2}-2 =-\frac{a}{\ell r^{n-1}}+\frac{r_{0}^{n}}{\ell r^{n}}+O(r^{-(n+1)}).
\end{align*}
Putting these into $\mathcal{E}$, we find that
\[
\mathcal{E}=-\frac{r_{0}^{n}}{\ell r^{n}}+O(r^{-2(n-1)}).
\]
Therefore
\begin{align*}
E(g) =\frac{1}{4\mbox{Vol}(\breve{T}^{n-1})}\lim_{r\to\infty}\int_{T_{r}^{n-1}}\mathcal{E}N\breve{e}^{2}\wedge\cdots\wedge\breve{e}^{n} =-\frac{\beta r_{0}^{n}}{8\pi\ell^{3}}=\frac{2G}{\lambda \ell} E_{{\tiny\mbox{HH}}}(g).
\end{align*}

In \cite{BCHMN}, Barzegar, Chru\'{s}ciel, H\"{o}rzinger, Maliborski and Nguyen proved the positive energy conjecture of Horowitz-Myers for the following metrics
\beQ
e^{2u} dr^{2}+e^{2v} d\phi^{2}+e^{2w}\sum_{i=1}^{n-2}(d\theta^{i})^{2}   \label{bchmn-g}
\eeQ
and $u$, $v$ and $w$ satisfy the following asymptotic conditions (Note $\ell=1$ in \cite{BCHMN}, here we choose arbitrary $\ell$)
\beQ
\begin{aligned}
u=-\ln \frac{r}{\ell} -\frac{1}{2}\ln \left(1-\frac{r_0 ^n}{r^n} \right)+\frac{u_n}{r^n} +o(r^{-n}),\quad
v=\ln \frac{r}{\ell} +\frac{1}{2}\ln \left(1-\frac{r_0 ^n}{r^n} \right)+\frac{v_n}{r^n} +o(r^{-n}),\quad
w=\ln r +\frac{w_n}{r^n}+o(r^{-n}).
\end{aligned}
\eeQ
Then we obtain
\beq
\begin{aligned}
e^{2u}=\frac{\ell ^2}{r^2}\left(1+\frac{r_0 ^n +2u_n}{r^n} +o(r^{-n})\right),\quad
e^{2v}=\frac{r^2}{\ell ^2}\left(1-\frac{r_0 ^n -2v_n}{r^n} +o(r^{-n})\right),\quad
e^{2w}=r^2 \left(1+\frac{2w_n}{r^n} +o(r^{-n})\right). \label{aij-bchmm}
\end{aligned}
\eeq

Compare (\ref{aij-g}) with (\ref{aij-bchmm}), we know that, if $a\neq 0$, the fall-offs of $g_{11}$ and $g_{22}$ of are slower than that of (\ref{bchmn-g}), thus we can not apply the main theorem for the positive energy conjecture of Horowitz-Myers in \cite{BCHMN} to $g$. In order to verify it for $g$, we require that $\phi$ has the same period in $g$ as the period $\frac{4\pi\ell^{2}}{n\overline{r}_{0}}$ in the following Horowitz-Myers metric
\[
g_{\tiny\mbox{HM}}=\frac{1}{\frac{r^{2}}{\ell^{2}}\left(1-\frac{\overline{r}_{0}^{n}}{r^{n}}\right)}dr^{2}+\frac{r^{2}}{\ell^{2}}\left(1-\frac{\overline{r}_{0}^{n}}{r^{n}}\right)d\phi^{2}+r^{2}\sum_{i=1}^{n-2}(d\theta^{i})^{2}.
\]
This implies
\begin{equation}
\overline{r}_{0}=\frac{r_{+}}{n}\left(n-1+\frac{r_{0}^{n}}{r_{+}^{n}}\right).\label{eq:HM-metric_r0}
\end{equation}

\begin{thm}
For metric $g$ of Horowitz-Myers type constructed in Theorem \ref{thm1} with fixed $r_0$, it holds that
\beQ
E(g)\geq E(g_{\tiny\mbox{HM}}),
\eeQ
with equality if and only if $a=0$, and $g$ is the Horowitz-Myers metric $g_{\tiny\mbox{HH}}$.
\end{thm}

\pf
For the Horowitz-Myers metric $g_{\tiny\mbox{HM}}$, we have
\[
E_{\tiny\mbox{HH}}(g_{\tiny\mbox{HM}})=-\frac{\lambda\overline{r}_{0}^{n-1}}{4Gn}<0.
\]
Combining this with (\ref{eq:HM-metric_r0}), it follows from (\ref{eq:HM-energy_2}) that
\[
E_{\tiny\mbox{HH}}(g)=-\frac{\lambda r_{0}^{n}}{4Gr_{+}\left(n-1+\frac{r_{0}^{n}}{r_{+}^{n}}\right)}=\left(\frac{n\frac{r_{0}}{r_{+}}}{n-1+\frac{r_{0}^{n}}{r_{+}^{n}}}\right)^{n}E_{\tiny\mbox{HH}}(g_{\tiny\mbox{HM}}).
\]
Since it is easy to see that $n-1+s^{n}\geq ns$ for $s\geq0$ and
the equality occurs if and only if $s=1$, applying this with $s=\frac{r_{0}}{r_{+}}$
to the above equation, we get
\[
E_{\tiny\mbox{HH}}(g)\geq E_{\tiny\mbox{HH}}(g_{\tiny\mbox{HM}}),
\]
with equality if and only if $r_{+}=r_{0}$. In this case, we have
$$0=V(r_{+})=V(r_{0})=a\frac{r_{0}^{n}}{\ell^{2}},$$ therefore $a=0$. By Proposition \ref{E-equal}, we complete the proof of the theorem.

\begin{acknowledgments}
The authors would like to thank referees for many valuable suggestions, especially for their comment on the conceptual understanding of the Horowitz-Myers conjecture in the early version of the paper.
The work is supported by Chinese NSF grants 11571345, 11701215, 11731001, the special foundation for Junwu and Guangxi Ba Gui Scholars, and HLM, NCMIS, CEMS, HCMS of Chinese Academy of Sciences.
\end{acknowledgments}

%\section*{DATA AVAILABILITY}
%Data sharing is not applicable to this article as no new data were created or analyzed in this study.

\bibliography{aipsamp}

\begin{thebibliography}{99}

\bibitem{ACD} M. Anderson, P. Chru\'{s}ciel, E. Delay, Non-trivial, static, geodesically complete, vacuum space-times with a negative cosmological constant. J. High Energy Phys. \textbf{10} (2002) 063.
\bibitem{BCHMN} H. Barzegar, P. Chru\'{s}ciel, M. H\"{o}rzinger, M. Maliborski, L. Nguyen, Remarks on the energy of asymptotically Horowitz-Myers metrics. Phys. Rev. D \textbf{101}, (2020) 024007.
\bibitem{CZ} J. Chen, X. Zhang, Metrics of Eguchi-Hanson types with the negative constant scalar curvature, arXiv:2007.15964 (2020).
\bibitem{D} D. Dold, Global dynamics of asymptotically locally AdS spacetimes with negative mass, Class. Quantum Grav. \textbf{35} (2018) 095012, 26 pp.
\bibitem{EH1} T. Eguchi, A.J. Hanson, Asymptotically flat self-dual solutions to Euclidean gravity, Phys. Lett. \textbf{74B} (1978) 249-251.
\bibitem{EH2} T. Eguchi, A.J. Hanson, Self-dual solutions to Euclidean gravity,  Ann. Phys. \textbf{120} (1979) 82-106.
\bibitem{GSW} G. Galloway, S. Surya, and E. Woolgar, On the Geometry and Mass of Static, Asymptotically AdS Spacetimes, and
the Uniqueness of the AdS Soliton, Commun. Math. Phys. \textbf{241} (2003) 1-25.
%\bibitem{HT} M. Henneaux, C. Teitelboim, Asymptotically anti-de Sitter spaces. Commun. Math. Phys. \textbf{98} (1985) 391-424.
\bibitem{HM} G. Horowitz, R. Myers,  AdS-CFT correspondence and a new positive energy conjecture for general relativity, Phys. Rev. \textbf{D 59}, no. 2 (1998): 026005.
\bibitem{L1} C. LeBrun, Counter-examples to the generalized positive action conjecture, Commun. Math. Phys. \textbf{118} (1988) 591-596.
\bibitem{L2} C. LeBrun, The Einstein-Maxwell Equations and Conformally K\"{a}hler Geometry, Commun. Math. Phys. \textbf{344} (2016) 621-653.
\bibitem{O} B. O\textquoteright Neill, Semi-Rimannian Geometry with Application to Relativity, vol. 103, Pure and Applied
Mathematics. Academic Press, New York (1983).
\bibitem{P} H. Pedersen, Eguchi-Hanson metrics with cosmological constant, Class. Quantum Grav. \textbf{2} (1985) 579-587.
\bibitem{Z} X. Zhang, Scalar flat metrics of Eguchi-Hanson type, Commun. Theor. Phys. (Beijing, China) \textbf{42} (2004) 235-237.
\end{thebibliography}

\end{document}